\documentclass[12pt,a4paper]{article} 
  
\usepackage{amsmath,amssymb} 
\usepackage{graphicx,color} 
\usepackage{url} 
  
\newcommand{\mbf}[1]{\protect\text{\boldmath$#1$}} 
\newcommand{\mbb}[1]{\mathbb{#1}} 
\newcommand{\cond}{\mathrm{cond}} 
 
\newcommand{\Tol}{\mathrm{Tol}\,} 
 
\newcommand{\ih}{\scalebox{0.67}[0.87]{$\Box$\hspace*{1pt}}} 
\newcommand{\USS}{\varXi_{\hspace{-0.5pt}uni}} 
\newcommand{\TSS}{\varXi_{\hspace{-0.5pt}tol}} 
\newcommand{\Ab}{(\mbf{A}, \mbf{b})} 
\newcommand{\IVE}{\mathrm{IVE}\,}
\newcommand{\un}[1]{\underline{#1}} 
\newcommand{\ov}[1]{\overline{#1}}

\newcommand{\m}{\mathrm{mid}\,} 
\renewcommand{\r}{\mathrm{rad}\,} 
  
\title{\bf A \ variability \ measure \ for \ estimates \\[2mm] 
          of parameters in interval data fitting} 
  
\author{%
     {\sc Sergey P. Shary}\\ 
     {\small Institute of Computational Technologies SB RAS}\\[-2pt] 
     {\small and Novosibirsk State University,}\\[-2pt] 
     {\small Novosibirk, Russia}\\[-2pt] 
     {\small E-mail: \tt shary@ict.nsc.ru} 
     } 
\date{} 
  
\begin{document}
  
\maketitle

\begin{abstract}
The paper presents a construction of a quantitative measure of variability for parameter
estimates in the data fitting problem under interval uncertainty. It shows the degree of 
variability and ambiguity of the estimate, and the need for its introduction is dictated 
by non-uniqueness of answers to the problems with interval data. A substantiation of
the new variability measure is given, its application and motivations are discussed.
Several examples and a series of numerical tests are considered, showing the features 
of the new characteristic and the specifics of its use. \\[2mm]
\hspace*{20.7mm} maximum compatibility method, strong compatibility, variability measure.\\[2mm]
\end{abstract}

\bigskip 
  
\section{Introduction and problem statement}

The purpose of this work is to present a quantitative variability measure for estimates 
of parameters of functional dependencies in the statistics of interval data. This is 
a relatively young branch of modern data science that does not rely on the probability 
theory, but makes extensive use of interval analysis methods (see, e.\,g., the surveys 
in \cite{ApplInteAnal,BoundApprHandbook,NguyenKreinWuXiang}). 
  
  
\begin{figure}[!htb]
\centering\small
\unitlength=1mm
\begin{picture}(84,70)
\put(0,-1){\includegraphics[width=84mm]{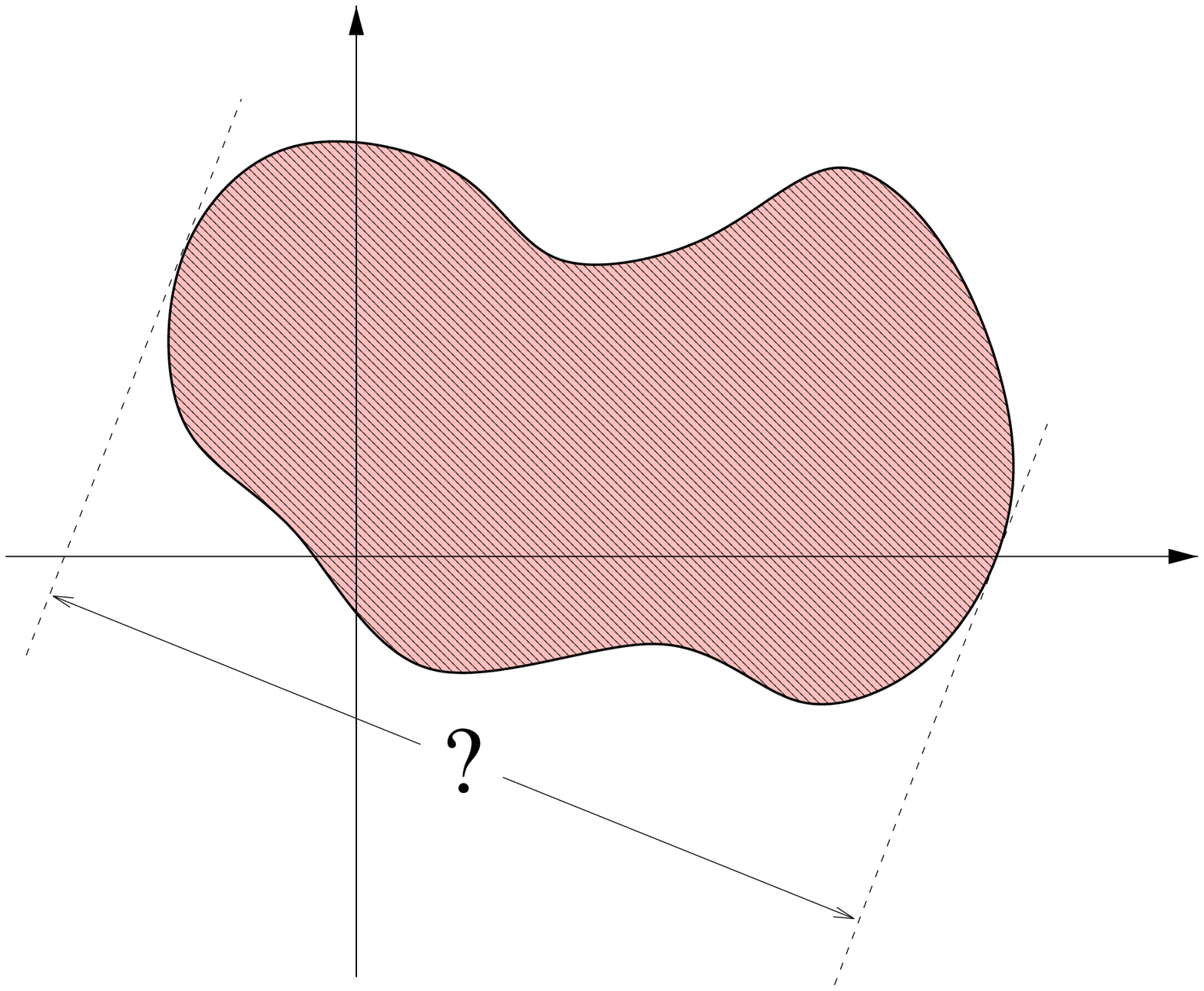}}
\end{picture}
\caption{\small A variability measure can be an estimate} 
of the size of the set of possible solutions. 
\label{SizeQuePic}
\end{figure}
  
  
By the term ``variability'', we understand the degree of variation and ambiguity 
of the estimate, and the need for its introduction is dictated by the fact that, 
in processing interval data, the answer is typically not unique. Usually, we get 
a whole set of different estimates that are equally consistent (compatible) with 
the source data  and, thus,  suitable  as solutions to the problem. The extent 
to which this set is large or small is, partly,  characterized by the term 
``variability''. In traditional probabilistic statistics, estimates of parameters 
are known to be random variables themselves, and the measure of their variability 
can be the variance of the estimates, mean absolute difference, median absolute 
deviation, average absolute deviation, and such like. What could be their 
analogues in the statistics of interval data? 
  
At first glance, the answer to this question seems quite obvious: it can be any value 
that characterizes the size of the set of solutions to the problem, if it is non-empty. 
We can even take an enclosure of the solution set obtained by an interval method. 
A certain disadvantage of this variant is the excessive detailing of the answer given 
as a box in $\mbb{R}^n$, a large amount of information that still needs to be 
``digested'' and reduced to a compact and expressive form. Sometimes, an interval 
estimate in the form of an axes-aligned box may inadequately represent the solution 
set. Another disadvantage is the complexity of finding such an estimate. 
  
It is desirable to have a relatively simple and efficiently computable quantity, 
expressed in a single number, because it would give a general aggregate view of 
the subject of interest. Similarly to variance and other probabilistic measures, 
it can serve as an approximate characteristic of the quality of parameter estimation. 
The greater the variability of an estimate, the less its certainty and the worse 
its quality, and this can serve as a basis for conclusions about the quality 
of the estimate. 
  
At the same time, the introduced variability measure should not be simply the ``size 
of the solution set''. If this solution set, for example, is unstable and changes 
abruptly with arbitrarily small changes in the data, then its size is, to some 
extent, misleading and disorienting (see example in Section 4). A practically useful 
variability measure should take into account this possible instability of the solution 
set to the problem and give us a robust value. 
  
  
\begin{figure}[!htb]
\centering\small
\unitlength=1mm
\begin{picture}(110,70)
\put(0,0){\includegraphics[width=110mm]{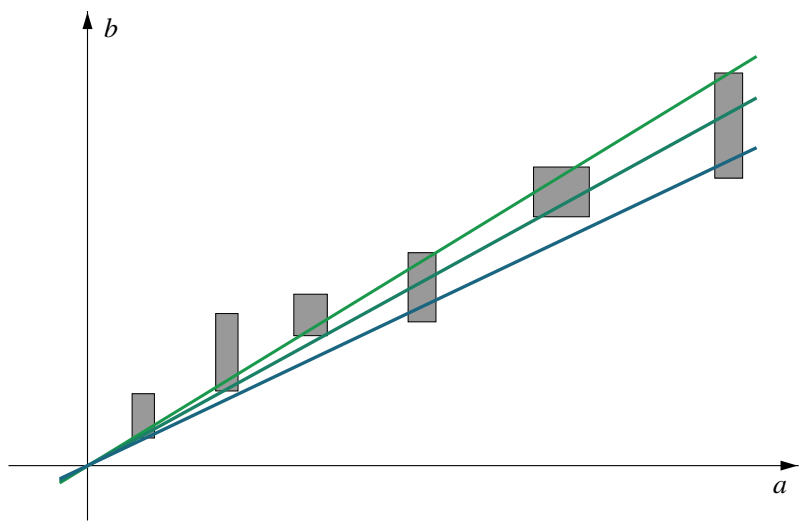}}
\end{picture}
\caption{\small An illustration for the data fitting problem under interval uncertainty.}
\label{IDataFitPic}
\end{figure}
  
  
In our article, we are within the framework of the data fitting problem (often 
called regression analysis problem): given results of measurements or observations, 
it is required to construct a functional dependence of a fixed type that ``best fit'' 
these data. Specifically, we need to determine the parameters $x_{1}$, $x_{2}$, 
\ldots, $x_{n}$ of a linear function of the form 
\begin{equation} 
\label{LinFunc} 
b = x_{1} a_{1} + \ldots + x_{n} a_{n} 
\end{equation} 
from a number of values of the independent variables $a_1$, $a_2$, \ldots, $a_n$ 
(also called \emph{exogenous}, \emph{explanatory}, \emph{predictor} or \emph{input} 
variables), and the corresponding values of the dependent variable $b$ (also called 
\emph{endogenous}, \emph{response}, \emph{criterion} or \emph{output} variable). 
Both $a_1$, $a_2$, \ldots, $a_n$ and $b$ are not known precisely, and we only have 
intervals of their possible values (see Fig.~\ref{IDataFitPic}). To find estimates 
of the coefficients $x_{1}$, $x_{2}$, \ldots, $x_{n}$, we use the so-called maximum 
compatibility method (previously called ``maximum consistency method''), which was 
proposed and developed in the works \cite{KreinovichShary,Shary-2012,Shary-2016, 
Shary-ADSAA} and others. After the estimates for $x_{1}$, $x_{2}$, \ldots, $x_{n}$ 
are found, we need to somehow evaluate their variability. Our article presents 
a construction of the variability measure in the above data fitting problem. 
  
Note that traditional methods of data fitting and regression analysis, such as 
the least squares method and its modifications, the least modulus method, etc., 
cannot be applied to the solution of our problem, since they are unsuitable for 
situations where the source data are intervals rather than points. 
  
  
\section{Formulation of the main results} 
  
  
\subsection{Maximum compatibility method and tolerable solution set}

The initial data for our problem is a set of values of independent and dependent 
variables for function \eqref{LinFunc}, which are obtained as a result of $m$ 
measurements (observations): 
\begin{equation} 
\label{EmpIData} 
\begin{array}{ccccc} 
\mbf{a}_{11}, & \mbf{a}_{12}, & \ldots & \mbf{a}_{1n}, & \mbf{b}_{1}, \\
\mbf{a}_{21}, & \mbf{a}_{22}, & \ldots & \mbf{a}_{2n}, & \mbf{b}_{2}, \\
  \vdots      &   \vdots      & \ddots &   \vdots      &  \vdots      \\
\mbf{a}_{m1}, & \mbf{a}_{m2}, & \ldots & \mbf{a}_{mn}, & \mbf{b}_{m}.
\end{array}
\end{equation} 
These are intervals as we assume that these data are inaccurate and have interval 
uncertainty due to measurement errors, etc. Both the data \eqref{EmpIData} and other 
interval values throughout the text are highlighted in bold mathematical font according 
to the informal international standard \cite{InteNotation}. The first index of the 
interval values from~\eqref{EmpIData} means the measurement number, and the second 
one, at $\mbf{a}_{ij}$'s, is the number of the independent variable that takes 
the corresponding value in this measurement. 
  
To find an estimate $(\hat{x}_{1}, \hat{x}_{2}, \ldots, \hat{x}_{n})$ of the parameters 
of the linear function \eqref{LinFunc}, we ``substitute'' data \eqref{EmpIData} into 
equality \eqref{LinFunc}, thus getting an interval system of linear algebraic equations 
\begin{equation} 
\label{InteLAS1} 
\arraycolsep=2pt 
\left\{ \ 
\begin{array}{ccccccccc} 
\mbf{a}_{11} x_1 &+& \mbf{a}_{12} x_2 &+& \ldots &+& \mbf{a}_{1n} x_n 
   &=& \mbf{b}_{1}, \\[1pt] 
\mbf{a}_{21} x_1 &+& \mbf{a}_{22} x_2 &+& \ldots &+& \mbf{a}_{2n} x_n 
   &=& \mbf{b}_{2}, \\[1pt] 
  \vdots         & &      \vdots      & & \ddots & &      \vdots 
   & &  \vdots      \\[1pt] 
\mbf{a}_{m1} x_1 &+& \mbf{a}_{m2} x_2 &+& \ldots &+& \mbf{a}_{mn} x_n 
   &=& \mbf{b}_{m}, 
\end{array} 
\right. 
\end{equation} 
or, briefly, 
\begin{equation}
\label{InteLAS2}
\mbf{A}x = \mbf{b}
\end{equation}
with an interval $m\times n$-matrix $\mbf{A} = (\mbf{a}_{ij})$ and interval 
$m$-vector $\mbf{b} = (\mbf{b}_{i})$ in the right-hand side. 
The sets of parameters which are compatible, in this or that sense, with 
the measurement data \eqref{EmpIData} form various solution sets for the equations 
system \eqref{InteLAS1}. The most popular of them are the \emph{united solution set} 
and \emph{tolerable solution set}. The united solution set, defined as 
\begin{equation*}
\USS(\mbf{A}, \mbf{b}) = \bigl\{\,x\in\mbb{R}^n
   \mid \text{ $Ax = b\,$ for some $A\in\mbf{A}$ and $b\in\mbf{b}$}\,\bigr\},
\end{equation*}
corresponds to the so-called weak compatibility between the parameters of function 
\eqref{LinFunc} and data \eqref{EmpIData} (see \cite{KreinovichShary,Shary-2012, 
Shary-2016}). The tolerable solution set, defined as 
\begin{equation*}
\TSS(\mbf{A}, \mbf{b}) = \bigl\{\,x\in\mbb{R}^n
   \mid \text{ $Ax\in\mbf{b}\,$ for each matrix $A\in\mbf{A}$}\,\bigr\}, 
\end{equation*}
corresponds to the so-called strong compatibility between the parameters of function 
\eqref{LinFunc} and data \eqref{EmpIData} (see \cite{Shary-ADSAA}). 
  
  
\begin{figure}[!htb]
\centering\small
\unitlength=1mm
\begin{picture}(110,65)
\put(0,0){\includegraphics[width=100mm]{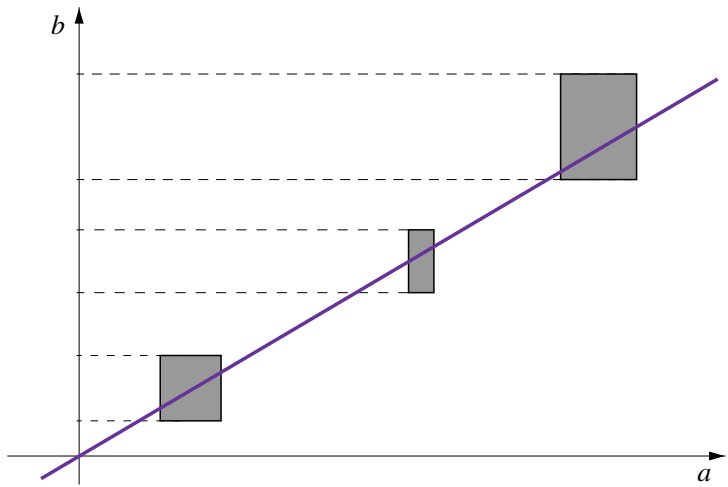}} 
\end{picture} 
\caption{\small An illustration of the strong compatibility between interval data 
         and a linear function.}
\label{StrCmpPic}
\end{figure}
  
  
Further, we assume that the solution to the data fitting problem for function 
\eqref{LinFunc} is found by the maximum compatibility method (see \cite{Shary-2012, 
Shary-2016, Shary-ADSAA}). As an estimate of the parameters of function \eqref{LinFunc}, 
it takes the maximum point of the \emph{recognizing functional}, a special function 
that gives a quantitative ``compatibility measure'' of this estimate with empirical 
data \eqref{EmpIData}. 
  
The maximum compatibility method has two versions, ``weak'' and ``strong'', that 
differ in understanding how exactly the interval data should be ``compatible'' with 
the function that we construct on them. Weak and strong compatibility reflect two 
different situations that may occur in data processing. In the weak version, it is 
required that the graph of the constructed function just intersects the measurement 
uncertainty boxes (see \cite{Shary-2012,Shary-2016}). The strong version implies 
more stringent condition: it requires that the function graph passes within the 
``corridors'' specified by the intervals $\mbf{b}_i$, $i = 1,2,\ldots,m$, for 
\emph{any} values of the independent variables $a_1$, $a_2$, \ldots, $a_n$ from 
the respective  intervals $\mbf{a}_{i1}$, $\mbf{a}_{i2}$, \ldots, $\mbf{a}_{in}$ 
obtained in the $i$-th measurement (see \cite{Shary-ADSAA}). This is illustrated 
in Fig.~\ref{StrCmpPic}, where the straight line of the function graph goes through 
the vertical faces of the measurement uncertainty boxes. The weak compatibility 
is shown in Fig.~\ref{IDataFitPic} by two upper straight lines. The lower line 
in Fig.~\ref{IDataFitPic} does not satisfy compatibility condition at all, 
neither weak nor strong, since it does not intersect some boxes. 
  
The ``strong version'' of the maximum compatibility method has a number of theoretical 
and practical advantages over the ``weak version''. These are polynomial complexity, 
robustness of estimates and their finite variability, the fact that the strong 
compatibility partially overcomes the so-called Demidenko paradox, etc. (see details 
in \cite{Shary-ADSAA}). Hence, we consider below a strong version of the maximum 
compatibility method, which corresponds to the tolerable solution set $\TSS(\mbf{A}, 
\mbf{b})$ for the interval system of equations \eqref{InteLAS2}. Its recognizing 
functional is usually denoted by ``Tol'', 
\begin{equation} 
\label{TolAbExpr} 
\Tol(x, \mbf{A}, \mbf{b})\;  = \,\min_{1\leq i\leq m}
  \left\{ \,\r\mbf{b}_i - \left|\; \m\mbf{b}_i - \sum_{j=1}^n
  \,\mbf{a}_{ij} x_{j} \,\right| \,\right\}, 
\end{equation} 
where 
\begin{equation*} 
\r\mbf{b}_{i} = \tfrac{1}{2}(\ov{\mbf{b}}_{i} - \un{\mbf{b}}_{i}), 
\hspace{23mm} 
\m\mbf{b}_{i} = \tfrac{1}{2}(\ov{\mbf{b}}_{i} + \un{\mbf{b}}_{i}) 
\end{equation*} 
are radii and midpoints of the components of the right-hand side $\mbf{b}$, 
the arithmetic operations inside the modulus in \eqref{TolAbExpr} are those 
of the classical interval arithmetic (see, e.\,g., \cite{ApplInteAnal, 
MooreBakerCloud,Neumaier}), and the modulus is understood as the maximum 
absolute value of the points from the interval, 
\begin{equation*} 
|\mbf{a}| = \max\,\{\,|a| \mid a\in\mbf{a}\,\} 
   = \max\,\bigl\{\,|\un{\mbf{a}}|, |\ov{\mbf{a}}|\,\bigr\}. 
\end{equation*} 
Typical graphs of the functional Tol for the one-dimensional case are shown 
in Fig.~\ref{MaxRecFuncPic} and Fig.~\ref{SteepFlatPic}. 
  
To solve the data fitting problem for the linear function \eqref{LinFunc} and data 
set \eqref{EmpIData}, it is necessary to find the unconstrained maximum, over all 
$x\in\mbb{R}^n$, of the functional $\Tol (x, \mbf{A}, \mbf{b})$, 
\begin{equation*}
\Tol(x, \mbf{A}, \mbf{b}) \rightarrow\max ,
\end{equation*} 
and the vector $\hat{x} = \arg\max_{x\in\mbb{R}^n}\,\Tol (x,\mbf{A}, \mbf{b})$ 
at which this maximum is attained provides an estimate of the parameters 
of function \eqref{LinFunc}. 
  
If $\max\,\Tol\geq 0$, then the solution set $\TSS(\mbf{A},\mbf{b})$, i.\,e., the set 
of parameters  strongly compatible  with the data  is non-empty,  and $\hat{x}\in\TSS 
(\mbf{A}, \mbf{b})$. If $\max\,\Tol < 0$, then the solution set $\TSS(\mbf{A}, \mbf{b})$ 
is empty and there do not exist parameters that are strongly compatible with data 
\eqref{EmpIData}. However, the argument $\hat{x}$ of $\max\,\Tol$ still provides 
the best compatibility of the constructed linear function with data \eqref{EmpIData} 
(more precisely, the least incompatibility). 
  
To conclude this subsection, we give a useful result on the tolerable solution set 
that allows us to investigate whether it is bounded or unbounded, i.\,e., whether 
the tolerable solution sets is finite in size or extends infinitely. 
  
\bigskip\noindent 
\textbf{Irene Sharaya's boundedness criterion} \cite{IreneRC} \ 
\textit{Let the tolerable solution set to an interval linear system $\mbf{A}x = \mbf{b}$ 
be nonempty. It is unbounded if and only if the matrix $\mbf{A}$ has linearly dependent 
noninterval columns.} 
  
\bigskip 
The criterion of boundedness shows that the tolerable solution set is unbounded, 
in fact, under exceptional circumstances, which are almost never fulfilled in practice, 
when working with actual interval data. That is, the tolerable solution set is mostly 
bounded, and the estimates obtained by the strong version of the maximum compatibility 
method almost always has finite variability.

  
\subsection{Variability measures}

As a quantity characterizing the variability of the estimate of the parameter vector 
$\hat{x} = (\hat{x}_{1}, \hat{x}_{2}, \ldots, \hat{x}_{n})$ in the linear function 
\eqref{LinFunc}, which is obtained by the maximum compatibility method 
from data \eqref{EmpIData}, we propose 
\begin{equation} 
\label{IVE} 
\IVE(\mbf{A}, \mbf{b})\; = \;\sqrt{n}\;\max_{\mbb{R}^n}\,\Tol \cdot 
   \Bigl(\;\min_{A\in\mbf{A}}\,\cond_{2}\,A\,\Bigr) \cdot \frac{\displaystyle 
   \bigl\|\,\arg\max_{\mbb{R}^n} \,\Tol\bigr\|_{2}}{\|\hat{\mbf{b}}\|_2}. 
\end{equation} 
In this formula, 
\begin{description} 
\item 
$n$ is the dimension of the parameter vector of function \eqref{LinFunc} under 
construction, 
\item 
$\|\cdot\|_2$ is the Euclidean norm (2-norm) of vectors from $\mbb{R}^n$, defined as 
\begin{equation*} 
\|x\|_2 \; = \; \left(\;\sum_{i=1}^n |x_{i}|^{2}\,\right)^{1/2}, 
\end{equation*} 
\item 
$\cond_{2}\,A$ is the spectral condition number of the matrix $A$, defined as 
\begin{equation*} 
\cond_{2}\,A\;  = \; \frac{\sigma_{\max}(A)}{\sigma_{\min}(A)}, 
\end{equation*} 
i.\,e., the ratio of the maximal $\sigma_{\max}(A)$ and minimal $\sigma_{\min}(A)$ 
singular values of $A$; it is an extension, to the rectangular case, of the concept 
of the condition number from computational linear algebra 
(see e.\,g. \cite{GolubVanLoan,Watkins}); 
\item 
$\hat{\mbf{b}}$ is a certain ``most representative'' point from the interval vector 
$\mbf{b}$, which is taken as 
\begin{equation} 
\label{hatbFormula}
\hat{\mbf{b}}\; = \;\tfrac{1}{2}(|\m\mbf{b} + \r\mbf{b}| + |\m\mbf{b} - \r\mbf{b}|), 
\end{equation} 
where the operations ``mid'' and ``rad'' are applied in componentwise manner. 
\end{description} 
  
  
\begin{figure}[htb] 
\centering\small 
\unitlength=1mm 
\begin{picture}(100,55) 
\put(0,0){\includegraphics[width=100mm]{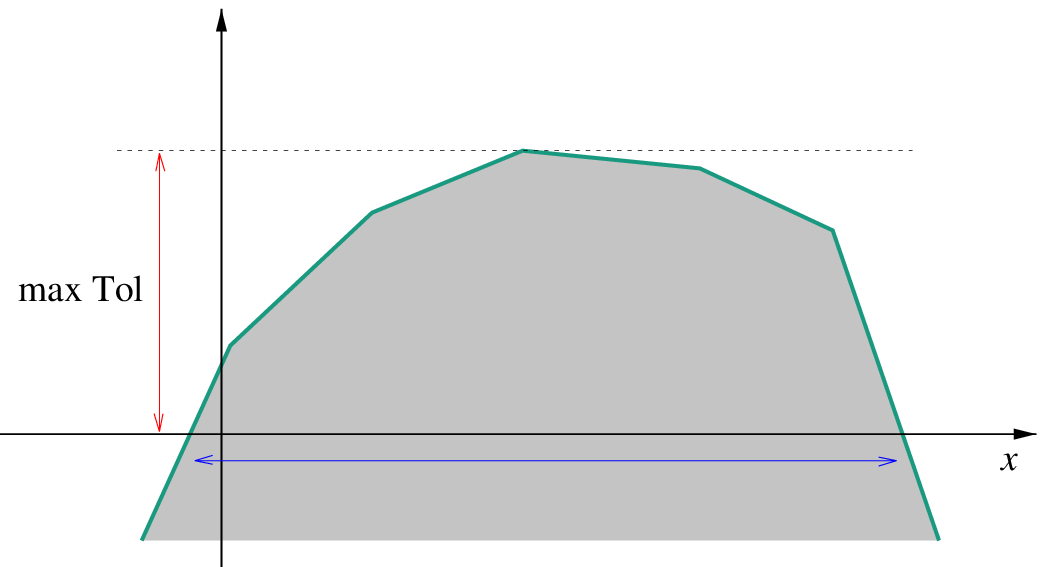}} 
\put(55,6){\large$\varXi_\mathit{tol}$} 
\end{picture}
\caption{\small The maximum value of the recognizing functional} 
\hspace*{0.5ex} gives an idea of the size of the tolerable solution set $\TSS$. 
\label{MaxRecFuncPic} 
\end{figure} 
  
  
Despite the definite formula \eqref{hatbFormula} for $\hat{\mbf{b}}$, it should be 
noted that the introduction of this point is, to a large extent, a matter of common 
sense. The general approach to the definition of $\hat{\mbf{b}}$ is that it must be 
a kind of ``most representative'' point from the right-hand side vector $\mbf{b}$, 
and in some situations this choice may be different from formula \eqref{hatbFormula}. 
For example, $\hat{\mbf{b}}$ can be a point result of the measurement, around which 
the uncertainty interval is built later, based on information about the accuracy 
of the measuring device. 
  
Apart from \eqref{IVE}, as a measure of relative variability of the parameter estimate, 
the value 
\begin{equation} 
\label{RVE} 
n\;\Bigl(\,\min_{A\in\mbf{A}}\,\cond_{2}A\,\Bigr) 
   \cdot\frac{\max_{\mbb{R}^n}\,\Tol}{\|\hat{\mbf{b}}\|_2}, 
\end{equation} 
can have a certain significance. Both IVE and value \eqref{RVE} are defined for 
interval linear systems \eqref{InteLAS2} with nonzero right-hand sides. They can 
take either positive real values or be infinite. The latter occurs in the only case of 
$\,\min_{A\in\mbf{A}}\,\cond_{2}A = \infty$, when all the point matrices $A\in\mbf{A}$ 
have incomplete rank, i.\,e., when $\sigma_{\min}(A) = 0$ for every $A\in\mbf{A}$. 
Then the variability measures are set to be infinite. 
  
The symbol IVE is built as an abbreviation of the phrase ``\un{i}nterval \un{v}ariability 
of the \un{e}stimate''. Below, we show that the value IVE adequately characterizes 
the size of non-empty tolerable solution set for a large class of practically important 
situations. But it is useful to discuss informal motivations that lead to the estimate 
IVE and to demonstrate that IVE has an intuitive, clear and even visual meaning. 
  
  
\begin{figure}[!htb]
\centering\small 
\unitlength=1mm
\begin{picture}(130,42)
\put(0,0){\includegraphics[width=130mm]{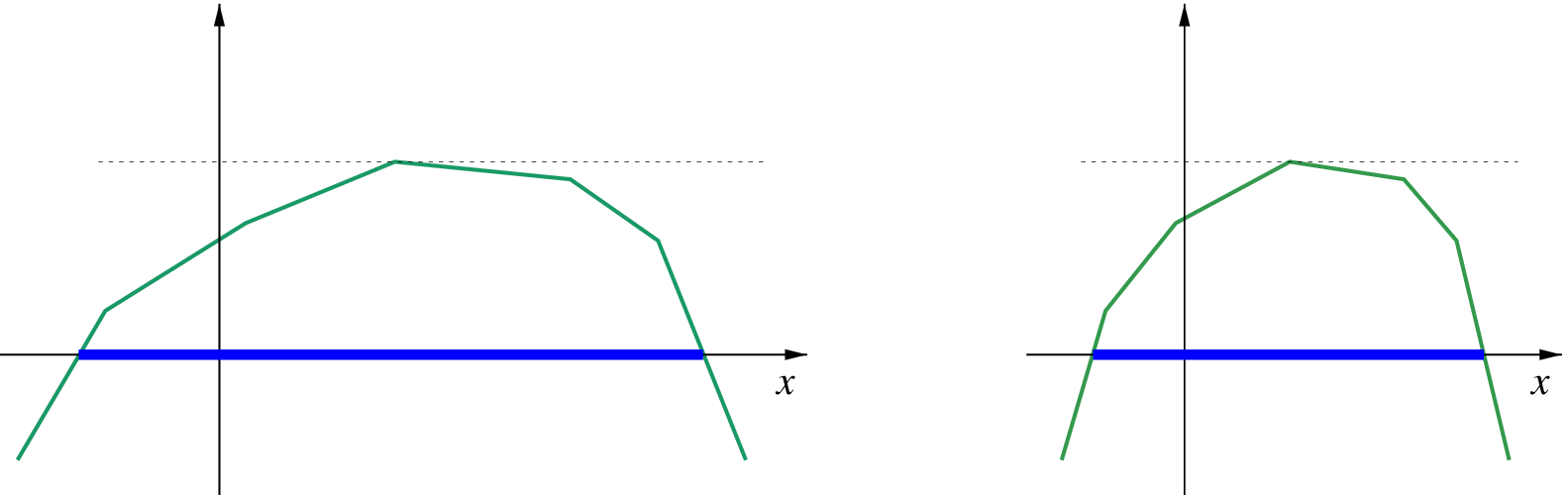}}
\put(35,7.5){\large$\varXi_\mathit{tol}$} 
\put(105,7.5){\large$\varXi_\mathit{tol}$} 
\end{picture}
\caption{\small In addition to the maximum of the recognizing functional, the size} 
\hspace*{0.8ex} of the tolerable solution set is also affected by ``steepness'' of the graph. 
\label{SteepFlatPic} 
\end{figure}
  
  
The tolerable solution set of an interval system of linear algebraic equations 
is the set of zero level of the recognizing functional Tol (see details in 
\cite{Shary-MCS}), or, in other words, the intersection of the hypograph of 
this functional with the coordinate plane $\Tol = 0$ (this is illustrated in 
Fig.~\ref{MaxRecFuncPic}). As a consequence, the magnitude of the maximum of 
the recognizing functional can, with other things being equal, be a measure of 
how extensive or narrow the tolerable solution set is. The more $\max\,\Tol$, 
the larger the size of the tolerable solution set, and vice versa. An additional 
factor that provides ``other things being equal'' is the slope (steepness) of 
pieces of hyperplanes of which the polyhedral graph of the functional Tol is 
compiled (these are straight lines in the 1D case in Fig.~\ref{MaxRecFuncPic} 
and Fig.~\ref{SteepFlatPic}). The slope of the hyperplanes is determined by 
the coefficients of the equations that define them, which are the endpoints 
of the data intervals \eqref{EmpIData}. The value of this slope is summarized 
in terms of the condition number of point matrices from the interval data matrix 
$\mbf{A}$. Finally, the multiplier 
\begin{equation*}
\frac{\|\arg\max\,\Tol\|_{2}}{\|\hat{\mbf{b}}\|_2} \ 
   = \  \frac{\|\hat{x}\|_2}{\|\hat{\mbf{b}}\|_2} 
\end{equation*}
is a scaling coefficient that helps to provide the commensurability of the final 
value with magnitudes of the solution, $\arg\max\,\Tol$, and the right-hand side 
vector of the equations system. Thus, formula \eqref{IVE} is obtained.

  
\section{A justification of the variability measure} 
\label{EstDerivSect}

Considering the most general case, we should assume that the number of measu\-rements 
$m$ may not coincide with the number $n$ of unknown parameters of the linear function 
\eqref{LinFunc}. In this section, we consider only the case $m\geq n$. In other words, 
the number of measurements (observations) made is not less than the number of function 
parameters. Then the interval system of linear equations \eqref{InteLAS2} is either 
square or tall (overdetermined). Of course, the data fitting problem makes sense 
for $m < n$ too, the maximum compatibility method also works for this case, and 
the variability measure IVE is then also applicable (see Section~4), but the latter 
still needs a separate substantiation.

  
\subsection{Estimates of perturbations of the solution \\
            to rectangular linear systems}

The starting point of our constructions justifying the choice of \eqref{IVE} exactly 
in the form described above is the well-known inequality that estimates perturbation 
$\Delta x$ of a nonzero solution $x$ to the system of linear algebraic equations 
$Ax = b$ depending on the change $\Delta b$ of the right-hand side $b$ (see, e.\,g., 
\cite{GolubVanLoan, Watkins}): 
\begin{equation} 
\label{CondIneq} 
\frac{\|\Delta x\|_2}{\|x\|_2} \ \leq \ 
   \cond_{2}\,A\,\cdot\frac{\|\Delta b\|_2}{\|b\|_2}. 
\end{equation} 
It is usually considered for square systems of linear equations, when $m = n$, 
but in the case of the Euclidean vector norm and the spectral condition number 
of matrices, this inequality holds true in the more general case with $m\geq n$. 
Naturally, estimate \eqref{CondIneq} makes sense only for $\sigma_{\min}(A)\neq 0$, 
when $\cond_{2} A < \infty$, i.\,e., when the matrix $A$ has full column rank. 
Let us briefly recall its derivation for this case. 
  
Given 
\begin{equation*}
Ax = b  \quad \text{ и }\quad  A(x + \Delta x) = b + \Delta b,
\end{equation*}
we have 
\begin{equation*}
A\Delta x = \Delta b.
\end{equation*}
Further, 
\begin{align*}
\displaystyle
\frac{\displaystyle\phantom{M}\frac{\|\Delta x\|_2}{\|x\|_2}\phantom{M}}%
     {\displaystyle\frac{\|\Delta b\|_2}{\|b\|_2}} \
&= \ \frac{\|\Delta x\|_{2}\,\|b\|_{2}\phantom{I}}%
     {\phantom{I}\|x\|_{2}\,\|\Delta b\|_{2}}
 = \ \frac{\|\Delta x\|_{2}\,\|Ax\|_{2}\phantom{I}}%
     {\phantom{I}\|x\|_{2}\,\|A\Delta x\|_{2}} \
 = \ \frac{\|\Delta x\|_2}{\|A\Delta x\|_2}\; \frac{\|Ax\|_2}{\|x\|_2} \\[3mm]
&\leq\;\max_{\Delta x\neq 0}\frac{\|\Delta x\|_2}{\|A\Delta x\|_2} \
     \max_{x\neq 0}\frac{\|Ax\|_2}{\|x\|_2} \
 = \ \left(\min_{\Delta x\neq 0}\frac{\|A\Delta x\|_2}{\|\Delta x\|_2}\right)^{-1} \
     \max_{x\neq 0}\frac{\|Ax\|_2}{\|x\|_2}                            \\[4mm]
&= \ \bigl(\sigma_{\min}(A)\bigr)^{-1} \,  \sigma_{\max}(A) \,
 = \ \cond_{2}(A)
\end{align*}
by virtue of the properties of the singular values (see e.\,g. \cite{HornJohnson,Watkins}). 
A comparison of the beginning and the end of this calculation leads to the inequality 
\eqref{CondIneq}, which, as is easy to understand, is attainable for some $x$ and $\Delta x$, 
or, equivalently, for some right-hand sides of $b$ and their perturbations $\Delta b$. 
Naturally, the above calculations and the resulting estimate make sense only 
for $\sigma_{\min}(A)\neq 0$.

  
\subsection{Interval systems with point matrices}

Let us consider an interval system of linear algebraic equations 
\begin{equation} 
\label{PointILAS} 
Ax = \mbf{b} 
\end{equation} 
with a point (noninterval) $m\times n$-matrix $A$, $m\geq n$, and an interval $m$-vector 
$\mbf{b}$ in the right-hand side. We assume that $A$ has full column rank and, therefore, 
$\cond_{2}\,A < \infty$. 
  
Suppose also that the tolerable solution set for system \eqref{PointILAS} is non-empty, 
i.\,e. $\TSS(A, \mbf{b}) = \bigl\{\,x\in\mbb{R}^n \mid Ax\in\mbf{b}\,\bigr\}\neq
\varnothing$. We need to quickly and with little effort estimate the size of this 
solution set, and our answer will be a ``radius type'' estimate for $\TSS(A, \mbf{b})$. 
More precisely, we are going to evaluate $\max \| x '- \hat{x}\|_2$ over all $x'\in\TSS 
(A,\mbf{b})$ and for a special fixed point $\hat{x}\in\TSS(A,\mbf{b})$, 
which is taken as 
\begin{equation*} 
\hat{x} \  = \  \arg\max_{x\in\mbb{R}^n}\,\Tol(x,A,\mbf{b}). 
\end{equation*} 
Recall that the argument $\hat{x}$ of the maximum of the recognizing functional 
for system \eqref{PointILAS} is an estimate of parameters of linear function 
\eqref{LinFunc} from empirical data. Strictly speaking, this point can be determined 
non-uniquely, but then let $\hat{x}$ be any one of the points at which the maximum 
is reached. 
  
Let $x'$ be a point in the tolerable solution set $\TSS(A,\mbf{b})$. How to evaluate 
$\|x' - \hat{x}\|_2$? It is clear that $x'$ and $\hat{x} $ are solutions of systems 
of linear algebraic equations with the matrix $A$ and some right-hand sides $b'$ and 
$\hat{b}$, respectively, from the interval vector $\mbf{b}$. If $\hat{x}\neq 0$ and 
$\hat{b}\neq 0$, then we can apply inequality \eqref{CondIneq}, considering 
a perturbation of the solution $\hat{x}$ to the system of linear algebraic equations 
$Ax = \hat{b}$. Then $\Delta x = x' - \hat{x}$, $\Delta b = b' - \hat{b}$, and we get 
\begin{equation*} 
\frac{\|x' - \hat{x}\|_2}{\|\hat{x}\|_2} \  \leq \; 
   \cond_{2}\,A \cdot\frac{\|b' - \hat{b}\|_2}{\|\hat{b}\|_2}, 
\end{equation*} 
from where the absolute estimate is obtained 
\begin{equation} 
\label{IntermEst} 
\|x' - \hat{x}\|_{2} \  \leq \; 
   \cond_{2}\,A \cdot\|\hat{x}\|_{2}\cdot \frac{\|b' - \hat{b}\|_2}{\|\hat{b}\|_2}. 
\end{equation} 
The point $\hat{x}$ is found as the result of maximization of the recognizing 
functional Tol, the point $\hat{b}$ coincides with $A\hat{x}$, the condition number 
$\cond_{2}\,A$ can be computed by well-developed standard procedures. Therefore, 
for practical work with inequality \eqref{IntermEst}, one need somehow evaluate 
$\|b' - \hat{b}\|_2$. 
  
But first, bearing in mind the further application of the deduced estimate in 
a situation where the matrix $A$ may vary, we somewhat roughen \eqref{IntermEst} 
by taking approximately $\|\hat{b}\|_{2} \approx \|\hat{\mbf{b}}\|_{2}$, that is, 
as the norm of the ``most representative'' point $\hat{\mbf{b}}$ of the interval 
vector $\mbf{b}$, which we defined in Section~2.2: 
\begin{equation*} 
\|\hat{b}\|_{2} \,\approx \, \|\hat{\mbf{b}}\|_2, 
\qquad \text{ where } \  \hat{\mbf{b}}\, = \,\tfrac{1}{2}\, 
   \bigl(\,|\m\mbf{b} + \r\mbf{b}| + |\m\mbf{b} - \r\mbf{b}|\,\bigr). 
\end{equation*} 
In doing this, some coarsening is allowed, so instead of \eqref{IntermEst} we write 
\begin{equation} 
\label{IntermEst2} 
\|x' - \hat{x}\|_{2} \  \lessapprox \; 
   \cond_{2}\,A \cdot\|\hat{x}\|_{2}\cdot \frac{\|\Delta b\|_2}{\|\hat{\mbf{b}}\|_2}. 
\end{equation} 
  
Now it is necessary to determine the increment of the right-hand side $\Delta b 
= b' - \hat{b}$. Its obvious upper bound is $2\,\r\mbf{b}$, but it is too crude. 
To get a more accurate estimate of $\Delta b$, we also consider, along with system 
\eqref{PointILAS}, a system of linear algebraic equations 
\begin{equation} 
\label{AuxILAS} 
Ax = \tilde{\mbf{b}}, 
\end{equation} 
for which the right-hand side is obtained by uniform ``compressing'' the interval 
vector $\mbf{b}$: 
\begin{equation} 
\label{ContractRHS}
\tilde{\mbf{b}}\, := \,\bigl[\,\un{\mbf{b}} +M, \ov{\mbf{b}} - M\,\bigr], 
\end{equation} 
where 
\begin{equation*} 
M \; := \; \max_{x\in\mbb{R}^n}\;\Tol(x, A, \mbf{b}) \  \geq \  0. 
\end{equation*} 
Since the maximum $M$ is reached for a certain value of the argument, $\hat{x}$, then 
\begin{equation*} 
M = \,\min_{1\leq i\leq m}
  \left\{ \,\r\mbf{b}_i - \left|\; \m\mbf{b}_i - \sum_{j=1}^n 
  \,\mbf{a}_{ij} \hat{x}_{j} \,\right| \,\right\}  \ 
  \leq \,\min_{1\leq i\leq m} \,\r\mbf{b}_i . 
\end{equation*} 
As a result, $\un{\mbf{b}} + M \leq \ov{\mbf{b}} - M$ in componentwise sense, and 
the endpoints in the interval vector \eqref{ContractRHS} do not ``overlap'' each other. 
  
But the properties of the recognizing functional imply that, for the interval 
system of linear algebraic equations \eqref{AuxILAS} with the right-hand side 
\eqref{ContractRHS}, the maximum of the recognizing functional is zero: 
\begin{equation*} 
\max_{x\in\mbb{R}^n}\;\Tol(x, A, \tilde{\mbf{b}}) \  = \  0. 
\end{equation*} 
Indeed, the values of $\r\mbf{b}_i$ are summands in all expressions in \eqref{TolAbExpr}, 
for which we take the minimum over $i = 1,2,\ldots,m$. Hence, if we simultaneously 
increase or decrease all $\r\mbf{b}_i$ by the same value, keeping the midpoints 
$\m\mbf{b}_i$ unchanged, then the total value of the recognizing functional will 
increase or decrease by exactly same value. In other words, if we take a constant 
$C\geq 0$ and the interval $m$-vector $\mbf{e} = ([-1, 1], \ldots, [-1, 1])^{\top}$, 
then, for the system $Ax = \mbf{b} + C\mbf{e}\,$ with all the right-hand sides expanded 
by $[-C, C]$, we have 
\begin{equation} 
\label{TolEvol} 
\Tol(x,A,\mbf{b} + C\mbf{e}) \  =  \  \Tol(x,A,\mbf{b}) + C. 
\end{equation} 
Therefore, 
\begin{equation} 
\label{MaxTolEvol} 
\max_{x\in\mbb{R}^n}\;\Tol(x,A,\mbf{b}+C\mbf{e}) \  = \  
   \max_{x\in\mbb{R}^n}\;\Tol(x,A,\mbf{b}) + C. 
\end{equation} 
The uniform narrowing of the right-hand side vector acts on the tolerable solution set 
and the recognizing functional in a completely similar way. If we narrow down all 
the components by the same value $M$, then the maximum of the recognizing functional 
of the new interval system also decreases by $M$. 
  
By virtue of the properties of the recognizing functional, the tolerable solution set 
$\TSS(A,\tilde{\mbf{b}})$ for system \eqref{AuxILAS} has empty interior (such sets are 
often called ``non-solid'' or ``meager''), which we will consider equivalent to ``having 
zero size''. Naturally, this is a simplifying assumption, since in reality the tolerable 
solution set corresponding to the zero maximum of the recognizing functional may be not 
a single-point set. But we still accept that. This implication is also supported by 
the fact that the situation with the zero maximum of the recognizing functional 
is unstable: the corresponding tolerable solution set can become empty 
with an arbitrarily small data perturbation (see Section~4). 
  
Another fact concerning the auxiliary system \eqref{AuxILAS} with the narrowed 
right-hand side, which follows from \eqref{TolEvol}--\eqref{MaxTolEvol}, is that 
the point $\hat{x}$ remains to be the argument of the maximum of the recognizing 
functional: 
\begin{equation*}
\hat{x} \  = \  \arg\max_{x\in\mbb{R}^n} \,\Tol (x,A,\tilde{\mbf{b}}). 
\end{equation*}
For this reason, the point $\hat{b} = A\hat{x}$ lies in the interval vector 
$\tilde{\mbf{b}}$ defined by \eqref{ContractRHS}. 
  
From what has been said, it follows that the solution set for the system $Ax = \mbf{b}$ 
is obtained from the solution set of the system $Ax = \tilde{\mbf{b}}$, which has 
``negligible size'' and for which $\max_{x\in\mbb{R}^n} \Tol(x,\mbf{A}, \tilde{\mbf{b}}) 
= 0$, through expanding the right-hand side vector $\tilde{\mbf{b}}$ in each component 
simultaneously by $[-M, M]$, where 
\begin{equation*}
M = \max_{x\in\mbb{R}^n}\;\Tol(x,\mbf{A}, \mbf{b}). 
\end{equation*} 
The interval vector $\tilde{\mbf{b}}\ni b$ may have non-zero size, but we put 
\begin{equation*} 
[-\Delta b, \Delta b] = ([-M, M], \ldots, [-M, M])^\top 
\end{equation*} 
in order to make our estimate 
\eqref{IntermEst2} attainable. Accordingly, in inequality \eqref{IntermEst2} we take 
\begin{equation*} 
\|\Delta b\| = \max_{x\in\mbb{R}^n}\;\Tol(x,\mbf{A},\mbf{b}), 
\end{equation*} 
if the Chebyshev norm ($\infty$-norm) is considered, or a value that differs from it 
by a corrective factor from the equivalence inequality for vector norms, if we take any 
other norm. As is known, for any vector $y\in\mbb{R}^n$ (see \cite{GolubVanLoan}) 
\begin{equation}
\label{EquivIneqs} 
\|y\|_{\infty}\leq \|y\|_{2} \leq \sqrt{n}\;\|y\|_{\infty}. 
\end{equation} 
Then 
\begin{equation} 
\label{PointEstim} 
\|x' - \hat{x}\|_{2} \  \lessapprox \, \sqrt{n} \  
   \,\cond_{2}A \cdot\bigl\|\,\arg\max_{x\in\mbb{R}^n}\,\Tol\bigr\|_{2} 
   \cdot \frac{\max_{x\in\mbb{R}^n}\Tol}{\|\hat{\mbf{b}}\|_2}. 
\end{equation} 
  
What happens if the matrix $A$ does not have a full column rank? Then, by virtue 
of the Irene Sharaya criterion, the nonempty tolerable solution set to the system 
\eqref{PointILAS} is unbounded. This is completely consistent with the fact that then 
$\cond_{2}A = \infty$ and the value of the variability measure IVE is infinite too.

  
\subsection{General interval systems}

Finally, we consider a general interval system of linear equations $\mbf{A}x = \mbf{b}$, 
with an essentially interval matrix, i.\,e., when $\r\mbf{A}\neq 0$.  In view of the 
properties of the tolerable solution set (see, e.\,g., \cite{Shary-MCS}), it can be 
represented as 
\begin{equation} 
\label{IxSSRepres} 
\TSS\Ab \  = \  \bigcap_{A\in\mbf{A}}\;\bigl\{\,x\in\mbb{R}^n \mid Ax\in\mbf{b}\,\bigr\} 
   \  =  \  \bigcap_{A\in\mbf{A}} \TSS(A, \mbf{b}), 
\end{equation} 
i.\,e., as the intersection of the solution sets to the individual systems $Ax = b$ 
with point matrices $A\in\mbf{A}$. 
  
For each interval linear system $Ax = \mbf{b}$ with $A\in\mbf{A}$, we have estimate 
\eqref{PointEstim}, if $A$ has full column rank. Otherwise, if the point matrix $A$ 
has incomplete column rank and the corresponding solution set $\TSS (A, \mbf{b})$ 
is unbounded, then we do not take it into account. Consequently, for the tolerable 
solution set of the system $\mbf{A}x = \mbf{b}$, which is the intersection of 
the solution sets $\TSS(A, \mbf{b})$ for all $A\in\mbf{A}$, the following should 
be true: 
\begin{equation} 
\label{RawIneq} 
\|x' - \hat{x}\|_{2} \  \lessapprox \  \min_{A\in\mbf{A}}\,
   \left\{\;\sqrt{n}\;\,\cond_{2}A \cdot \bigl\|\,\arg\max_{x\in\mbb{R}^n}\,\Tol\bigr\|_{2} 
   \cdot\frac{\max_{x\in\mbb{R}^n}\Tol}{\|\hat{\mbf{b}}\|_2}\,\right\}. 
\end{equation} 
The transition from representation \eqref{IxSSRepres} to inequality \eqref{RawIneq} can 
be both very accurate and rather crude (as can be seen from considering the intersection 
of two 1D intervals). It all depends on the size of the intersection of the solution sets 
of individual systems $Ax = \mbf{b}$. On the other hand, the amount of this intersection 
is indirectly characterized by the magnitude of $\max_{x\in\mbb{R}^n}\Tol$. 
  
Taking the above facts into account, we perform approximate estimation of the right-hand 
side of inequality \eqref{RawIneq} by moving the minimum over $A\in\mbf{A}$ through 
the curly brackets. First of all, we evaluate the factor $\|\arg\max_{x\in\mbb{R}^n}\, 
\Tol\|_{2}$, which changes to the smallest extent, by the constant available to us after 
the numerical solution of the data fitting problem: 
\begin{equation} 
\bigl\|\arg\max_{x\in\mbb{R}^n}\,\Tol(x, A, \mbf{b})\bigr\|_{2}\, 
    \approx \; \mathrm{const}\; = \;\bigl\|\arg\max_{x\in\mbb{R}^n}\, 
    \Tol(x,\mbf{A}, \mbf{b})\bigr\|_{2}.  \label{ArgMaxConst} 
\end{equation} 
Next, the minimum of $\cond_{2}A$ naturally turns to $\min\cond_{2}A$, and the most important 
factor $\max_{x\in\mbb{R}^n}\,\Tol(x, A, \mbf{b})$ will be changed to $\max_{x\in\mbb{R}^n} 
\Tol(x,\mbf{A}, \mbf{b})$. This choice (as well as \eqref{ArgMaxConst}) is rather rigidly 
determined by the following reasons. 
  
The expression for our variability measure should 
preserve its simplicity and be uniform for all cases and situations. In particular, 
if the interval matrix $\mbf{A}$ squeezes to a point matrix $A$, then our measure should 
turn to the estimate \eqref{PointEstim} for the point case. Finally, if $\max\,\Tol = 0$, 
then our measure must be zero too, since the size of the (stable) tolerable solution set 
is also zero, and our variability measure should reliably detect such situations. 
All this taken together leads to the estimate 
\begin{equation} 
\label{IntCondIneq} 
\|x' - \hat{x}\|_{2} \;\, \lessapprox \  \sqrt{n} \  \Bigl(\;\min_{A\in\mbf{A}} 
   \,\cond_{2}A \,\Bigr) \cdot \bigl\|\,\arg\max_{x\in\mbb{R}^n}\,\Tol\bigr\|_{2} 
   \cdot\frac{\max_{x\in\mbb{R}^n}\Tol}{\|\hat{\mbf{b}}\|_2}. 
\end{equation} 
  
The same estimate as \eqref{IntCondIneq}, by virtue of the equivalence inequality 
\eqref{EquivIneqs}, is also true for the Chebyshev norm: 
\begin{equation*} 
\max_{x'\in\TSS\Ab} \|x' - \hat{x}\|_{\infty} \ 
\lessapprox \  \sqrt{n} \;\, \Bigl(\;\min_{A\in\mbf{A}}\,\cond_{2}\,A\,\Bigr) 
  \cdot \bigl\|\,\arg\max_{x\in\mbb{R}^n}\Tol\bigr\|_{2} 
  \cdot \frac{\max_{x\in\mbb{R}^n}\,\Tol}{\|\hat{\mbf{b}}\|_2}. 
\end{equation*} 
This completes the rationale for \eqref{IVE}. 
  
If we want to evaluate the relative size of the tolerable solution set, expressing 
it in ratio to the norm of its points, then it is reasonable to take $\hat{x} = \arg 
\max_{x\in\mbb{R}^n}\Tol$ as the ``most typical'' point from the tolerable solution 
set $\TSS\Ab$. Using \eqref{EquivIneqs} again, we obtain 
\begin{equation*} 
\frac{\max_{x'\in\TSS\Ab} \|x' - x''\|_{\infty}}{\|\hat{x}\|_{\infty}} \ 
   \lessapprox \  n\; \Bigl(\;\min_{A\in\mbf{A}}\,\cond_{2}A\,\Bigr) 
   \cdot \frac{\max_{x\in\mbb{R}^n}\,\Tol}{\|\hat{\mbf{b}}\|_2}. 
\end{equation*} 
This gives value \eqref{RVE}.

  
\section{Numerical examples and some tests}

\begin{figure}[htb]
\centering\small
\unitlength=1mm
\begin{picture}(70,46)
\put(0,0){\includegraphics[width=70mm]{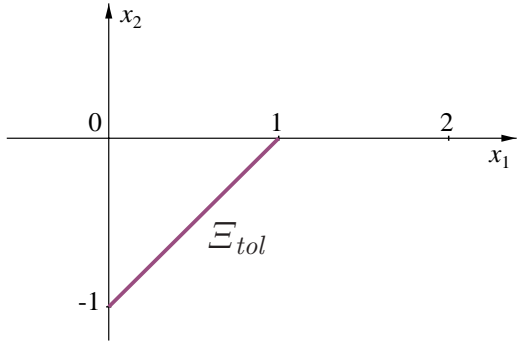}} 
\put(28,12){\large$\varXi_\mathit{tol}$} 
\end{picture}
\caption{\small  The tolerable solution set 
         for the interval equations systems \eqref{UnstableISys}.} 
\label{Exmp0Pic} 
\end{figure}

First of all, we consider an example of unstable tolerable solution set that changes 
abruptly with small perturbations in the system of equations. For all interval 
$2\times 2$-systems of linear algebraic equations of the form 
\begin{equation} 
\label{UnstableISys} 
\begin{pmatrix}
[-1, 1] & [-1, 1] \\[2pt] 
   1    &   -1    \\[2pt] 
\end{pmatrix}
\begin{pmatrix}
x_{1} \\[2pt] x_{2}
\end{pmatrix}
=
\begin{pmatrix}
{[-1, 1]} \\[2pt]
{[1, 1 + \eta]} 
\end{pmatrix}, \qquad \eta \geq 0,  
\end{equation}
the tolerable solution sets are the same: this is the straight line segment joining 
the points $(0,-1)$ and $(1,0)$ and depicted in Fig.~\ref{Exmp0Pic}. The diameter of 
the solution set is essentially non-zero (namely, $\sqrt{2}$), but the unconstrained 
maximum of the recognizing functional Tol for all such systems is zero, and it is 
attained at the point $(0.5,-0.5)$. 
  
  
At the same time, any arbitrarily small increase in the lower endpoint of the interval 
$[1,1+\eta]$ in the right-hand side of the second equation makes the tolerable solution 
set empty. An arbitrarily small reduction of the upper endpoint of the interval $[-1, 1]$, 
located in the first component of the right-hand side vector, produces a similar effect. 
It turns out that the maximum value of the recognizing functional Tol characterizes 
very precisely the instability of the original solution set and the zero size of 
the solution sets of perturbed systems. 
  
As the second example, we consider the problem of constructing a linear function 
of two variables $a_1$ and $a_2$, 
\begin{equation}
\label{2LinFunc}
b = x_{1} a_{1} + x_{2} a_{2}, 
\end{equation}
from the interval data obtained in 3 measurements: 
\begin{equation*}
\arraycolsep=4mm
\begin{array}{c|ccc}
 & \mbf{a}_{1} & \mbf{a}_{2} & \mbf{b} \\
\hline \\[-3mm]
1 & [98, 100] & [99, 101] & [190, 210] \\[3pt]
2 & [97, 99]  & [98, 100] & [200, 220] \\[3pt]
3 & [96, 98]  & [97, 99]  & [190, 210]
\end{array}
\end{equation*} 
Note that in these data the three-dimensional uncertainty boxes of measurements~1 
and 2, as well as 2 and 3, substantially ``overlap'' each other: their intersections 
are boxes with non-empty interiors, the sizes of which are comparable to the sizes 
of the original data boxes. 
  
  
\begin{figure}[htb]
\centering\small
\unitlength=1mm
\begin{picture}(96,78)
\put(0,0){\includegraphics[width=96mm]{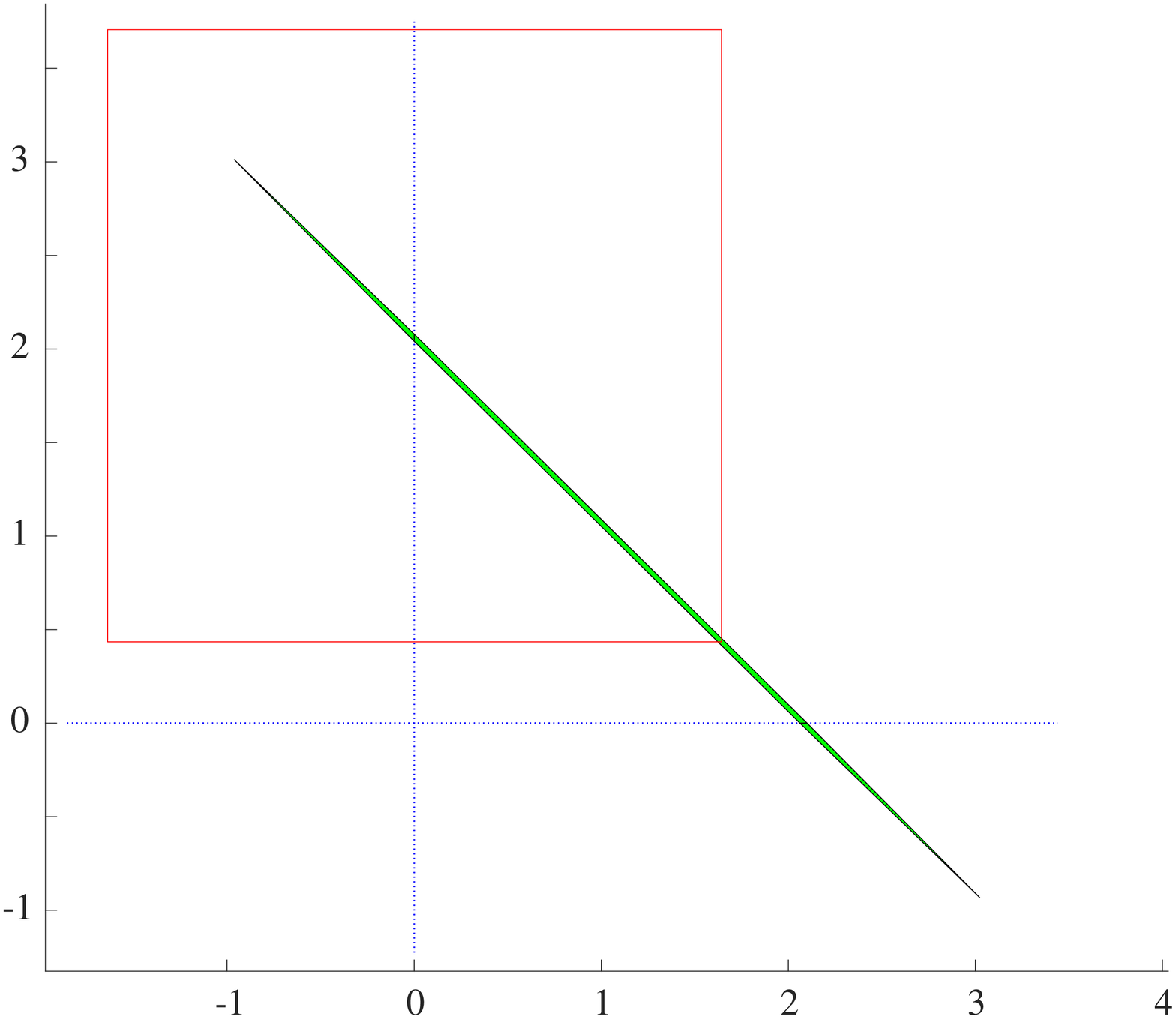}}
\put(55,1){\mbox{$x_1$}} \put(3,36){\mbox{$x_2$}} 
\put(60,70){\mbox{$\hat{\mbf{x}}$}}
\end{picture}
\caption{\small\,The\, tolerable\, solution\, set\, of\, the\, system\, 
                of\, equations\, \eqref{1stSampleISys}} 
in comparison with the box constructed by using the estimate IVE. 
\label{1stExmpPic} 
\end{figure}
  
  
To determine the coefficients $x_1$ and $x_2$, we compose an interval 
$3\times 2$-system of linear algebraic equations 
\begin{equation} 
\label{1stSampleISys} 
\begin{pmatrix}
[98, 100] & [99, 101] \\[2pt] 
[97, 99]  & [98, 100] \\[2pt] 
[96, 98]  & [97, 99] 
\end{pmatrix}
\begin{pmatrix}
x_{1} \\[2pt] x_{2}
\end{pmatrix}
=
\begin{pmatrix}
{[190, 210]} \\[2pt]
{[200, 220]} \\[2pt]
{[190, 210]}
\end{pmatrix}.
\end{equation}
Its matrix has incomplete rank, since its member is a point matrix with the rank 1:
\begin{equation} 
\label{Rank1Matr} 
\begin{pmatrix}
98 & 99 \\
98 & 99 \\
98 & 99
\end{pmatrix}.
\end{equation} 
The united solution set for system \eqref{1stSampleISys} is unbounded, therefore it is 
hardly possible to determine, with certainty, the coefficients of the linear function 
\eqref{2LinFunc} satisfying the weak compatibility between parameters and data 
(see Section~2). However, the interval matrix of system \eqref{1stSampleISys} does not 
contain linearly dependent point columns, and therefore, according to the Irene Sharaya 
criterion \cite{IreneRC} (see Section~2.1), the tolerable solution set is bounded. 
It is depicted in Fig.~\ref{1stExmpPic}, which is drawn by the procedure \texttt{EqnTol2D} 
from the package \texttt{IntLinInc2D} \cite{IreneVisu}. The minimum spectral condition 
number of the point matrices contained in the interval matrix of \eqref{1stSampleISys} 
is $103.83$, and it is reached on the matrix 
\begin{equation*} 
\begin{pmatrix}
100 &  99 \\
 97 & 100 \\
 96 &  99
\end{pmatrix}. 
\end{equation*} 
This result can be obtained, for example, using the simulated annealing algorithm 
on the set of point matrices contained in the interval matrix of \eqref{1stSampleISys}. 
  
Numerical solution of the maximization problem for the recognizing functional Tol can 
be carried out within MATLAB environment, using the free program \texttt{tolsolvty2.m} 
(available from the author's page at ResearchGate \cite{TOLSOLVTY2}). It implements 
a modified version of the so-called $r$-algorithms for non-differentiable optimization 
proposed and developed by N.Z.\,Shor and N.G.\,Zhurbenko \cite{ShorZhurbenko}. Using 
the precision settings specified in it ``by default'', we get $\max\,\Tol = 1.9095$, 
which is reached at $\hat{x} = (5.1857\cdot 10^{-7}, 2.0603)^\top$. Then, 
\begin{equation*}
\IVE = \sqrt{2}\cdot 1.9095\cdot 103.83\cdot 
   \frac{\|\hat{x}\|_2}{\sqrt{200^2 + 210^2 + 200^2}} = 1.6399. 
\end{equation*}
  
Interval hull of the tolerable solution set for system \eqref{1stSampleISys} 
(that is, its optimal interval enclosure) is the box 
\begin{equation*}
\begin{pmatrix}
[-0.9620, 3.0227] \\[2pt]
[-0.9320, 3.0127]
\end{pmatrix},
\end{equation*} 
which can also be found by the procedure \texttt{EqnTol2D}. We see that the value 
of IVE is in satisfactory agreement with the radii of the components of the optimal 
estimate of the solution set, equal to $1.9924$ and $1.9724$ respectively. 
  
In the maximum compatibility method, the argument $\hat{x} = \arg\max_{x\in\mbb{R}^n}
\Tol$ of the unconstrained maximum of the recognizing functional plays a crucial role, 
and, in fact, our variability estimate IVE relies heavily on it. This is why it makes 
sense to look at the box $\hat{\mbf{x}}$ with the components $[\hat{x}_{i} - \IVE, 
\hat{x}_{i} + \IVE]$, $i = 1,2$. It is also depicted in Fig.~\ref{1stExmpPic}, and 
the substantial asymmetry of its location relative to the solution set is, of course, 
explained by the specific position of the center, the point $\hat{x}$, as well as 
the ill-conditioning of the point systems from \eqref{1stSampleISys}. With other data, 
the box $\hat{\mbf{x}}$ estimates the tolerable solution sets significantly better
(see further). 
  
Next, we give an example of the opposite type (in a sense, dual to the previous 
example), where a linear function of three variables 
\begin{equation} 
\label{3vLinFun} 
b = x_{1} a_{1} + x_{2} a_{2} + x_{2} a_{3} 
\end{equation} 
is to be constructed from the data of two experiments summarized below: 
\begin{equation*} 
\arraycolsep=4mm 
\begin{array}{c|cccc} 
  & \mbf{a}_1 & \mbf{a}_2 &  \mbf{a}_3 & \mbf{b}   \\ 
\hline \\[-3mm] 
1 & [98, 100] & [97,  99] &  [96, 98] & [190, 210] \\[3pt] 
2 & [99, 101] & [98, 100] &  [97, 99] & [200, 220] 
\end{array} 
\end{equation*} 
To find the parameters of function \eqref{3vLinFun}, we come to an underdetermined 
interval system of linear algebraic equations 
\begin{equation}
\label{2ndSampleISys} 
\begin{pmatrix}
[98, 100] & [97, 99]  & [96, 98] \\[2pt]
[99, 101] & [98, 100] & [97, 99]
\end{pmatrix}
\begin{pmatrix}
x_{1} \\[1pt] x_{2} \\[1pt] x_{3}
\end{pmatrix}
=
\begin{pmatrix}
{[190, 210]} \\[2pt]
{[200, 220]}
\end{pmatrix}. 
\end{equation} 
Its matrix is the transposed matrix of system \eqref{1stSampleISys}, so 
$\min_{A\in\mbf{A}}\cond_{2}\,A$ is the same for it. Also, the matrix of system 
\eqref{2ndSampleISys} contains a point matrix of the incomplete rank~1, which is 
transposed for \eqref{Rank1Matr} (and many more such matrices). 
  
  
\begin{figure}[htb]
\centering\small
\unitlength=1mm
\begin{picture}(93,82)
\put(0,0){\includegraphics[width=93mm]{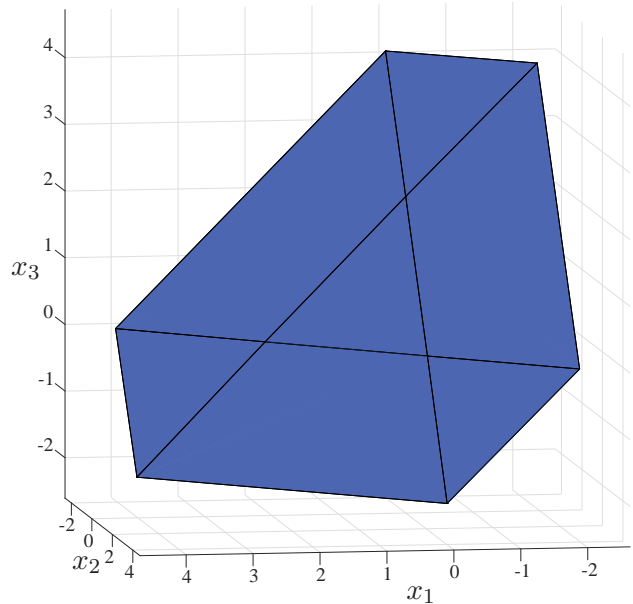}}
\put(55,2){\mbox{$x_1$}} \put(11,6){\mbox{$x_2$}}
\put(3,45){\mbox{$x_3$}}
\end{picture}
\caption{\small The tolerable solution set 
         for the interval equations system \eqref{2ndSampleISys}.} 
\label{2ndExmpPic} 
\end{figure}
  
  
Again, the united solution set for system \eqref{2ndSampleISys} is unbounded, and 
it is difficult (if at all possible) to determine the coefficients of the linear 
function \eqref{3vLinFun}, relying on the weak compatibility between parameters and 
data, due to ``infinite variability'' of the resulting estimate. Nevertheless, in 
these adverse conditions, the nonempty tolerable solution set to the interval system 
of equations \eqref{2ndSampleISys} is bounded by virtue of the Irene Sharaya criterion 
\cite{IreneRC} (see Section~2.1). In Fig.~\ref{2ndExmpPic}, the tolerable solution set 
is depicted as a thin hexagonal plate. Computation of the maximum of the recognizing 
functional for this system using the code \texttt{tolsolvty2.m} gives the value 
$\max \Tol = 3.9698$, which is reached at the point 
\begin{equation*}
\hat{x} = \arg\max \Tol  
   =  \,\bigl(\,2.0603, 3\cdot 10^{-6}, 2.1\cdot 10^{-6}\,\bigr)^{\top}. 
\end{equation*} 
It can be taken as an estimate of the coefficients in \eqref{3vLinFun}. Then 
the varibility measure of the above estimate is 
\begin{equation*} 
\IVE = \sqrt{2}\cdot 3.9698\cdot 103.83\cdot 
   \frac{\|\hat{x}\|_2}{\sqrt{200^2 + 210^2}} = 4.1413. 
\end{equation*} 
  
Interval hull (optimal interval enclosure) of the tolerable solution set for system 
\eqref{2ndSampleISys} is the box 
\begin{equation*}
\begin{pmatrix}
[-1.9747, 4.0302] \\[2pt]
[-1.9899, 4.0759] \\[2pt]
[-1.9949, 4.1071]
\end{pmatrix},
\end{equation*}
which can also be computed by the procedure \texttt{EqnTolR3}. The radii of 
the components of this interval vector are $3.0024$, $3.0329$, $3.0510$ respectively, 
which is also not very different from the value of IVE. The example shows that 
the value IVE works even in the case of $m < n$, when the number of measurements 
is less than the number of parameters to be determined. But a rigorous justification 
of this fact is waiting for its study. 
  
To conclude the section, we present, in Table~1, the results of numerical tests 
for the interval linear $n\times n$-system 
\begin{equation} 
\label{NeumNelSys} 
\left( 
\begin{array}{cccc} 
 \theta  & {[0, 2]} & \cdots & {[0, 2]} \\[1pt] 
{[0, 2]} &  \theta  & \cdots & {[0, 2]} \\[1pt] 
 \vdots  &  \vdots  & \ddots &  \vdots  \\[1pt] 
{[0, 2]} & {[0, 2]} & \cdots &  \theta 
\end{array} 
\right)
\; 
\left( 
\begin{array}{@{\,}c@{\,}}
x_{1}  \\[1pt] 
x_{2}  \\[1pt] 
\vdots \\[1pt]
x_{n} 
\end{array} 
\right) 
= 
\left(
\begin{array}{@{\;}c@{\;}}
{[1, K]} \\[1pt] 
{[1, K]} \\[1pt] 
 \vdots  \\[1pt]
{[1, K]} 
\end{array} 
\right), 
\end{equation} 
with various $n$ and $K$. System \eqref{NeumNelSys} resembles that proposed 
in \cite{Neumaier}, having exactly same matrix. But the right-hand sides were taken 
as positive intervals $[1, K]$, since the original balanced intervals $[-1, 1]$ 
in the system from \cite{Neumaier} make the tolerable solution set ``too symmetric''. 
  
  
\begin{table}[h!] 
\caption{Results of the computational tests with system \eqref{NeumNelSys}}
\small\centering 
\begin{tabular}{|c|c|c|c||c|c|c|c|}
\hline 
$\theta\rule[-3mm]{0mm}{8mm}$ & $\IVE$ & $\|\,\r(\,\ih\TSS)\|_\infty$ & 
   $\|\,\r(\,\ih\TSS)\|_2$ & 
$\theta\rule[-3mm]{0mm}{8mm}$ & $\IVE$ & $\|\,\r(\,\ih\TSS)\|_\infty$ & 
   $\|\,\r(\,\ih\TSS)\|_2$ \\ 
\hline 
\multicolumn{4}{|c|}{$n = 5$, $K = 10$} & \multicolumn{4}{|c|}{$n = 10$, $K = 10$} \\ 
2.0   &  1.019   &  1.25  &  2.795 &  6.0   &  0.894  &  0.5   &  1.581 \\ 
4.0   &  1.081   &  0.875 &  1.957 &  9.0   &  1.491  &  0.389 &  1.230 \\ 
6.0   &  0.786   &  0.639 &  1.429 &  12.0  &  0.582  &  0.313 &  0.988 \\ 
8.0   &  0.681   &  0.5   &  1.118 &  15.0  &  0.495  &  0.26  &  0.822 \\ 
10.0  &  0.534   &  0.41  &  0.917 &  20.0  &  0.396  &  0.203 &  0.640 \\ 
\hline 
\multicolumn{4}{|c|}{$n = 5$, $K = 20$} & \multicolumn{4}{|c|}{$n = 10$, $K = 20$} \\ 
2.0   &  2.953   &  3.75  &  8.385 &  6.0   & 2.489   &  1.333 &  4.216 \\ 
4.0   &  2.698   &  2.125 &  4.752 &  9.0   & 1.831   &  0.944 &  2.987 \\ 
6.0   &  2.015   &  1.472 &  3.292 &  12.0  & 1.432   &  0.729 &  2.306 \\ 
8.0   &  1.591   &  1.125 &  2.516 &  15.0  & 1.255   &  0.593 &  1.876 \\ 
10.0  &  1.378   &  0.91  &  2.035 &  20.0  & 0.985   &  0.453 &  1.431 \\ 
\hline 
\end{tabular} 
\end{table} 
  
  
The interval matrix of system \eqref{NeumNelSys} is known to be regular if and only 
if $\theta > n$ for even $n$ and $\theta > \sqrt{n^2 - 1}$ for odd $n$ \cite{Neumaier}. 
Consequently, in Table~1, the first two rows that correspond to each separate case 
of $n$ and $K$ refer to systems with singular matrices. As the parameter $\theta$ grows, 
the matrix of the system becomes regular and better conditioned. 
  
The values of IVE are compared with $\|\,\r(\,\ih\TSS)\|_\infty$ and 
$\|\,\r(\,\ih\TSS)\|_2$, that is, with the Chebyshev norm (max-norm) and 
the Euclidean norm of the radius of the interval hull of the tolerable solution set 
(denoted as $\ih\TSS$). We can see that, with the exception of two cases corresponding 
to $n = 5$ and $K = 10, 20$, the values of IVE are always within the two-sided bounds 
given by $\|\,\r(\,\ih\TSS)\|_\infty$ (lower bound) and $\|\,\r(\,\ih\TSS)\|_2$ 
(upper bound). And that is reasonable. Overall, our numerical experiments confirm 
the adequacy of the new measure of variability, which gives quite satisfactory 
approximate estimates of the size of the tolerable solution sets in various 
situations.

  
\section{Discussion}

IVE is the first measure of variability proposed in the statistics of interval data,
for estimation using the maximum compatibility method, and we can not compare IVE 
with similar other measures, since they simply do not exist. However, it is useful 
to correlate the estimate IVE with the ideal mathematical characteristics of the 
solution set, such as its diameter, in terms of computational convenience and 
laboriousness. 
  
First of all, IVE reflects instabilities in the solution set better than the diameter 
(see the first example in Section~4). An instability of the tolerable solution set 
for an interval linear system arises in the case when the maximum value of the 
recognizing functional is zero, $\max\Tol = 0$. Then the tolerable solution set 
can be either a single-point or an extended set with non-zero diameter and empty 
interior \cite{Shary-MCS}. After an arbitrarily small perturbation of data, the latter 
situation can abruptly turn into the empty solution set. In any case, this phenomenon 
is signaled by the zero value of the maximum of the recognizing functional. 
The corresponding variability measure IVE is also zero, which is quite natural: 
it makes sense to show only ``stable size'' of the solution set. The equality of IVE 
to zero or ``almost  zero'' thus allows us to diagnose unstable cases. 
  
Next, the problem of computing the diameter, in 2-norm, of the tolerable solution 
set to an interval linear system of equations is NP-hard in general. This follows 
from its reducibility to the quadratic programming problem with linear constraints 
(see \cite{Vavasis}). Indeed, the membership of a point in the tolerable solution set 
to an interval $m\times n$-system of equations is determined by a system of linear 
inequalities of the size $2m\times 2n$ (the Rohn theorem \cite{Rohn}). To compute 
the diameter of the tolerable solution set in 2-norm, we have to maximize the quadratic 
objective function $\|x' - x''\|^2_2$ over all pairs of points $x'$, $x''$ from 
the tolerable solution set, i.\,e. satisfying $2m\times 2n$-systems of linear 
inequalities. So, computing the diameter of the tolerable solution set is not easy. 
  
The diameter of the interval hull of the tolerable solution set can be computed more 
simply, but it is not better than IVE in any case, since the interval hull is not 
the solution set itself, and the coarsening resulted from such a replacement may 
be large. 
  
Calculation of IVE by formula \eqref{IVE} requires solving the data fitting problem, 
that is, finding $\max\,\Tol$ and $\arg\max\,\Tol$, and then we need to evaluate 
the minimum of the condition number of the matrices from the interval data matrix. 
In turn, the recognizing functional Tol is a concave piecewise linear function 
\cite{Shary-MCS}, so computing its maximum is polynomially complex. The author 
efficiently solves it by nonsmooth optimization methods developed in recent decades, 
in particular, using $r$-algorithms proposed by N.Z.\,Shor \cite{ShorZhurbenko}, or 
using separating plane algorithms (see, e.\,g., \cite{Nurminski, Vorontsova}). 
The most difficult part in calculating IVE is thus evaluating the minimum condition 
number of point matrices from a given interval matrix. 
  
Computing the exact minimum of the condition number is not simple, but to address 
practical problems  which will apply  the value IVE,  it  is  sufficient  to have 
an approximate estimate for $\min_{A\in\mbf{A}}\,\cond_{2}\,A$ from above. This 
follows from our considerations in Section~3.3. Sometimes, it is not necessary 
to compute $\min\,\cond_{2}\,A$ at all, if we have to compare, with each other, 
the variability of the estimates obtained for the same data matrix $\mbf{A}$. 
  
If the interval matrix is ``sufficiently narrow'', being not very different from 
a point matrix, then we can approximate 
\begin{equation} 
\label{MinApprxMid} 
\min_{A\in\mbf{A}}\,\cond_{2}\,A \;\approx \;\cond_{2}(\m\mbf{A}). 
\end{equation} 
But in general, this recipe may work poorly, since the left and right sides 
of the approximate equality \eqref{MinApprxMid} can be quite different. 
In the examples with systems \eqref{1stSampleISys} and \eqref{2ndSampleISys} 
from Section~4, the condition number of the midpoint matrix is $2.38\cdot 10^{4}$, 
and, using the simplified formula \eqref{MinApprxMid}, we are mistaken in evaluating 
the variability measure IVE by more than 20 times. 
  
In the general case, for a more accurate evaluation of $\min\,\cond_{2}\,A$, we can 
use popular evolutionary optimization methods, such as a genetic algorithm, simulated 
annealing, particle swarm optimization, etc., within the interval matrix $\mbf{A}$. 
In the numerical experiments from Section 4, the minimum of the condition number 
was found using the standard program of simulated annealing from free computer math 
system Scilab. 
  
Note that there is a fundamental difference between computing the variability measure 
IVE and computing the diameter of the tolerable solution set: in the first case, 
we calculate a minimum, while in the second we have to find a maximum. 
Using traditional optimization methods and various heuristics, in the first case 
we compute an approximation to the minimum from above, and in the second case we find 
an approximation to the maximum from below. If we want to get, with our variability 
measure, a guaranteed outer estimation of the solution set, then the upper estimate, 
which is obtained by calculating the minimum in IVE, is more preferable. 
  
There exists one more viewpoint at the variability measure IVE. 
  
In traditional probabilistic statistics, the phenomenon of collinearity of data 
(also called ``multicollinearity'') plays a large role. It is the presence of a linear 
dependence between the input (predictor) variables of the regression model. The $k$ 
variables of the model in question are usually called \emph{collinear} if the vectors 
representing them lie in a~linear space of dimension less than $k$ \cite{RegrDiagnos}, 
so that one of these vectors is a linear combination of the others. In practice, such 
exact collinearity of data is rare, but real computational problems in data fitting 
often starts when the data vectors are ``almost linearly dependent''. Then the parameter 
estimates are unstable, which leads to increased statistical uncertainty, i.\,e., 
an increase in the variance of the estimates. 
  
According to modern views, the collinearity of data is most adequately described 
by the condition number of the matrix made up of these data (see, e.\,g., 
\cite{RegrDiagnos}, Chapter~3). In this sense, our IVE is, in fact, a measure 
of the collinearity of the data, corrected with the help of the actual value of 
the estimate and compatibility of this estimate with the data (which is indicated 
by the maximal value of the recognizing functional). The minimum over all $\cond_{2}A$ 
for $A\in\mbf{A}$ is taken due to the specifics of the strong compatibility of 
parameters and data, and it agrees well with the regularizing role of 
the tolerable solution set (see~\cite{Shary-arXiv}). 
  
With this interpretation, IVE makes sense even with a negative maximum of 
the recognizing functional, max Tol, when the tolerable solution set is empty 
and the parameters of the linear function \eqref{LinFunc}, which are strongly 
compatible with the data, do not exist. The absolute value of IVE still shows, 
up to a certain scaling coefficient, a measure of the collinearity of the data 
(a measure of their ill-conditioning), and the negative sign indicates the status 
of the solution to the problem, i.\,e., that the parameter vector computed is not 
strongly compatible with the data, but only provides the best possible approximation 
for the input data of the problem. 
  
The immediate goal of further research is to justify the use of IVE for underdetermined 
situations, where the number $m$ of observations is less than the number $n$ of parameters 
to be determined. The maximum compatibility method works well in this case too, we get 
parameter estimates and we can calculate their values of IVE, but its application needs 
to be justified, at least at the same level of rigor as was done in this work for 
$m\geq n$. 
  
  
\paragraph{Aknowledgements} 
The author is grateful to Alexander Bazhenov, Sergey Kumkov, and Sergei Zhilin for 
stimulating discussions and valuable comments on the work. Also, the author thanks 
the anonymous reviewers for their constructive criticism and good suggestions.


  
\end{document}